\documentclass[reqno,11pt]{amsart}
\usepackage{amssymb, amsmath}
\usepackage{graphicx}
\usepackage{cite}
\usepackage{tensor}
\usepackage{appendix}
\usepackage{comment}

\usepackage{subfig, tikz}
\usetikzlibrary{decorations.pathmorphing}
\usepackage{graphics}

\usepackage{bm}
\usepackage{enumerate}
\theoremstyle{plain}

\theoremstyle{definition}

\def\be{\begin{equation}}
\def\ee{\end{equation}}
\def\bea{\begin{eqnarray}}
\def\eea{\end{eqnarray}}

\newcounter{mnotecount}[section]

\renewcommand{\themnotecount}{\thesection.\arabic{mnotecount}}

\newcommand{\mnote}[1]
{\protect{\stepcounter{mnotecount}}$^{\mbox{\footnotesize
$
\bullet$\themnotecount}}$ \marginpar{
\raggedright\tiny\em
$\!\!\!\!\!\!\,\bullet$\themnotecount: #1} }

\numberwithin{equation}{section}

\begin{document}
\title{Conformal Geodesics Cannot Spiral - Erratum}
\author{Peter Cameron}
\address{Department of Mathematics\\
Imperial College London\\
London SW7 2AZ, UK.
}
\email{p.cameron24@imperial.ac.uk}
\author{Maciej Dunajski}
\address{Department of Applied Mathematics and Theoretical Physics\\ 
University of Cambridge\\ Wilberforce Road, Cambridge CB3 0WA, UK.\\and\\
Faculty of Physics, University of Warsaw Pasteura 5, 02-093 Warsaw, Poland}
\email{m.dunajski@damtp.cam.ac.uk}
\author{Paul Tod}
\address{The Mathematical Institute\\
Oxford University\\
Woodstock Road, Oxford OX2 6GG\\ UK.
}
\email{tod@maths.ox.ac.uk}
\date{29 January 2026}
\maketitle
\begin{abstract}
Wojciech Kami\'nski  has provided a  non real--analytic counterexample to our claim in \cite{CDT}
that conformal geodesics cannot spiral. This erratum illustrates how the proof of Lemma 4.6 \cite{CDT} (on which our claim was based) fails.
\end{abstract}
Wojciech Kami\'nski \cite{Kaminski} has provided a 
non real--analytic counterexample to our claim in \cite{CDT}
that conformal geodesics cannot spiral (Theorem 2.3). 
Apart from the missing assumption of real--analyticity the error in this paper is contained in Lemma 4.6.
In the proof of this Lemma it is asserted that the bound on the function $F$ is a continuous function of $(p,u_0,\hat{a}_0)$. However, Kami\'nski's example shows that this is false and and consequently the size $R(p,u_0)$ (Definition 4.1) of the heart $H_{p,u_0}$ (Equation (3.21)), is not bounded away from 0 on any neighbourhood of the point towards which the conformal geodesic spirals. In Kami\'nski's example $R(\mu(t),\dot{\mu}(t))\rightarrow0$ along the spiralling conformal geodesic $\mu$, where $t$ is a parameter such that $g(\dot{\mu}(t),\dot{\mu}(t))=1$.

  We now illustrate how the proof of Lemma 4.6 fails. Let $p\in M$ and $u_0$ be a unit tangent vector at $p$. Suppose there exists an $\hat{a}_0$ with $g(u_0,\hat{a}_0)=0$ such that the conformal geodesic $\gamma$ with initial conditions $(p,u_0,\hat{a}_0)$ first exits the heart $H_{p,u_0}$ at $p$ (Figure 1). Now suppose there exists a sequence $(\gamma_n)_{n=1}^\infty$ of conformal geodesics with initial conditions $(p_n,u_n,\hat{a}_n)$ such that (writing $u_{0,n}=u_{0,n}^i\partial_i$ in local co-ordinates)
\begin{itemize}
\item $p_n\rightarrow p$ and $u_{0, n}^{i}\rightarrow u_0^i$ as $n\rightarrow\infty$
    \item the point at which $\gamma_n$ first exits the heart $H_{p_n,u_n}$ is not $p_n$.
\end{itemize}
It follows from the first point above that, for a fixed range of arc length parameter, the conformal geodesics $\gamma_n$ become uniformly close to $\gamma$. However, since $\gamma_n$ no longer first leaves the heart $H_{p_n,u_{0,n}}$ at $p_n$ it is possible that $R(p_n,u_{0,n})\nrightarrow R(p,u_0)$ (Figure 2). 
 Other results the paper, in particular the definition of the exponential map (Definition 3.1) as well as Theorem 3.4 are unaffected by the above error.

\begin{center}
\includegraphics[scale=0.20]{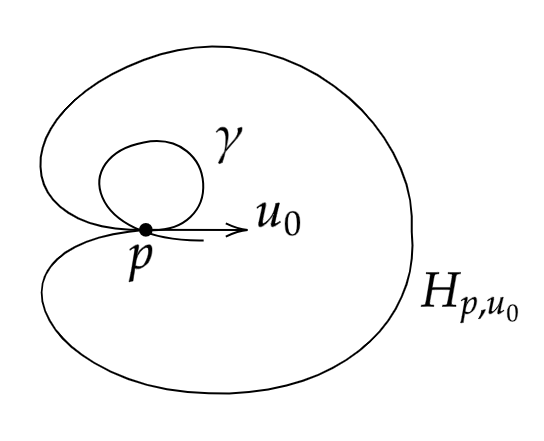}\quad  \includegraphics[scale=0.20]{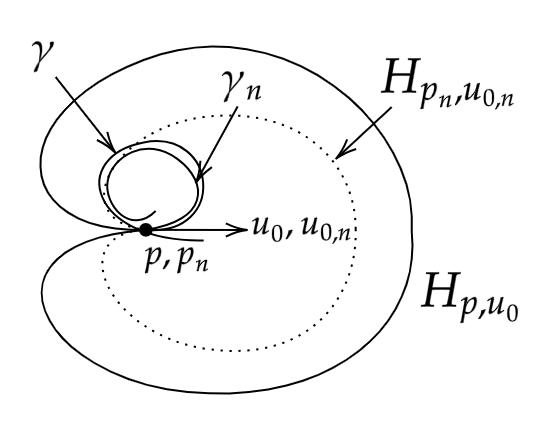}
\begin{center}
{\em Figure 1.\;\;\;\;\;\;\;\;\;\;\;\;\;\;\;\;\;\;\;\;\;\;\;\;\;\; Figure 2.}
\end{center}
\end{center}
{\em Figure 1}: The self-intersecting conformal geodesic $\gamma$ with initial unit tangent vector $u_0$ at $p$ first exits the heart $H_{p,u_0}$ at $p$.\\
{\em Figure 2}:
This figure shows {\em Figure 1} with a conformal geodesic $\gamma_n$ overlaid. As $n\rightarrow\infty$, $\gamma_n$ becomes uniformly close to $\gamma$ (for a fixed arc length parameter range) but $R(p_n,u_{0,n})\nrightarrow R(p,u_{0})$. For illustrative purposes we have drawn $p$ and $p_n$ as well as $u_0$ and $u_{0,n}$ to co-incide in this Figure. Formally this should be interpreted as $\gamma$ and $H_{p,u_0}$ plotted with respect to co-ordinates centred at $p$, while $\gamma_n$ and $H_{p_n,u_{0,n}}$ are plotted with respect to different co-ordinates centred at $p_n$.
\subsubsection*{Acknowledgements} We are grateful to Wojciech Kami\'nski for informing us of his example, and to  Mike Eastwood for discussions and advice. 

\end{document}